\documentclass[amsbook,12pt]{article}
\usepackage{epsfig, amsfonts,amssymb, amsmath}
\usepackage[usenames]{color}
\usepackage{picinpar}
\newcommand{\nl}{\mbox{}\\}

\begin{document}
%
%
%
\mbox{} \vspace{-1.500cm} \\
\begin{center}
{\Large \bf %
An inequality for solutions of the} \\
\mbox{} \vspace{-0.300cm} \\
\mbox{\Large \bf %
Navier-Stokes equations in $\mathbb{R}^{n}$} \\
\nl
\mbox{} \vspace{-0.300cm} \\
{\large \sc T. Hagstrom,}$\mbox{}^{\!\!\:\!1}$
{\large \sc J. Lorenz,}$\mbox{}^{\!\!\;\!2}$
{\large \sc J. P. Zingano}$\mbox{}^{\:\!3}$
{\large \sc and P. R. Zingano}$\mbox{}^{\:\!3}$ \\
\mbox{} \vspace{-0.200cm} \\
$\mbox{}^{1}${\small
Department of Mathematics} \\
\mbox{} \vspace{-0.685cm} \\
{\small
Southern Methodist University} \\
\mbox{} \vspace{-0.685cm} \\
{\small
Dallas, TX 75275-0235, USA} \\
\mbox{} \vspace{-0.390cm} \\
$\mbox{}^{2}${\small
Department of Mathematics and Statistics} \\
\mbox{} \vspace{-0.685cm} \\
{\small
University of New Mexico} \\
\mbox{} \vspace{-0.685cm} \\
{\small
Albuquerque, NM 87131-0001, USA} \\
\mbox{} \vspace{-0.390cm} \\
$\mbox{}^{3}${\small Departamento de Matem\'atica Pura e Aplicada} \\
\mbox{} \vspace{-0.670cm} \\
{\small Universidade Federal do Rio Grande do Sul} \\
\mbox{} \vspace{-0.670cm} \\
{\small Porto Alegre, RS 91509-900, Brazil} \\
\nl
\mbox{} \vspace{-0.500cm} \\
{\bf Abstract} \\
\mbox{} \vspace{-0.525cm} \\
%
%
\begin{minipage}[t]{11.350cm}
{\small
\mbox{} \hspace{+0.200cm}
We obtain a new inequality
that holds for general Leray solutions
of the incompressible Navier-Stokes equations
in $ \mathbb{R}^{n} \!$ ($n \leq 4$).
This recovers
important results
previously obtained
by other authors
regarding the time
decay of solution derivatives
(of arbitrary order). \\
}
\end{minipage}
\end{center}
%
%
%
%
\nl
\mbox{} \vspace{-0.775cm} \\
\mbox{} \hspace{+0.500cm}
{\bf 2010 AMS Subject Classification:}
{\small 35B40} (primary),
{\small 35D30}, {\small 35Q30}, {\small 76D05} \\
\mbox{} \vspace{-0.250cm} \\
\mbox{} \hspace{+0.500cm}
{\bf Keywords:}
time decay of solution derivatives,
Leray global weak solutions, \\
\mbox{} \hspace{+0.985cm}
incompressible Navier-Stokes equations,
%
%
%
Schonbek-Wiegner estimates \\
\nl
\setcounter{page}{1}
\mbox{} \vspace{-0.550cm} \\
%
%
%

%
%

{\bf 1. Introduction} \\
\mbox{} \vspace{-0.750cm} \\

In this note we derive a fundamental
new inequality for
general Leray solutions
of the incompressible Navier-Stokes equations
(in dimension $ n \leq 4 $),
that is,
global solutions
$ {\displaystyle
\;\!
\mbox{\boldmath $ u $}(\cdot,t)
\in
L^{\infty}((0,\infty), \:\!\mbox{\boldmath $L$}^{2}_{\sigma}(\mathbb{R}^{n}))
\,\cap
} $
$ {\displaystyle
L^{2}((0,\infty), \:\!\dot{\mbox{\boldmath $H$}}\mbox{}^{\!\;\!1}\!\;\!(\mathbb{R}^{n}))
\cap \;\!
C_{\mbox{\scriptsize w}}
([\;\!0, \infty), \mbox{\boldmath $L$}^{2}(\mathbb{R}^{n}))
} $
%
of the fluid flow system \\
\mbox{} \vspace{-0.725cm} \\
\begin{equation}
\tag{1.1$a$}
\mbox{\boldmath $u$}_{t}
\;\!+\,
\mbox{\boldmath $u$} \!\cdot\! \nabla \mbox{\boldmath $u$}
\,+\,
\nabla p
\;=\;
\nu \,
\Delta\:\! \mbox{\boldmath $u$},
%
%
\end{equation}
\mbox{} \vspace{-0.975cm} \\
\begin{equation}
\tag{1.1$b$}
\nabla \!\cdot  \mbox{\boldmath $u$}(\cdot,t) \,=\, 0,
\end{equation}
\mbox{} \vspace{-0.975cm} \\
\begin{equation}
\tag{1.1$c$}
\mbox{\boldmath $u$}(\cdot,0)
\,=\,
\mbox{\boldmath $u$}_0 \in
\mbox{\boldmath $L$}^{2}_{\sigma}(\mathbb{R}^{n}),
\end{equation}
\mbox{} \vspace{-0.175cm} \\
that satisfy
the generalized energy inequality 

\mbox{} \vspace{-1.250cm} \\
\begin{equation}
\tag{1.2}
\|\, \mbox{\boldmath $u$}(\cdot,t) \,
\|_{\mbox{}_{\scriptstyle L^{2}(\mathbb{R}^{n})}}
  ^{\:\!2}
\!\:\!+\,
2 \, \nu \!\!\;\!
\int_{\mbox{}_{\scriptstyle \!\:\!s}}^{\;\!t}
\!\!\:\!
\|\, D \,\!\mbox{\boldmath $u$}(\cdot,\tau) \,
\|_{\mbox{}_{\scriptstyle L^{2}(\mathbb{R}^{n})}}
  ^{\:\!2}
d\tau
\,\leq\:
\|\, \mbox{\boldmath $u$}(\cdot,s) \,
\|_{\mbox{}_{\scriptstyle L^{2}(\mathbb{R}^{n})}}
  ^{\:\!2}
\!\:\!,
\quad
\forall \;\;\! t \geq s
\end{equation}
\mbox{} \vspace{-0.150cm} \\
for a.e.$\;\!\;\!s \geq 0 $,
including $ s = 0 $.
Such solutions were first constructed
by Leray \cite{Leray1933, Leray1934} \linebreak
for $ n \leq 3 $,
and later by other authors
with different methods
and more general space dimension,
see e.g.\;\cite{FujitaKato1964, 
Kato1984, Ladyzhenskaya1969, 
ShinbrotKaniel1966, Sohr2001, Temam1984}.
%
In (1.1) above,
$ \nu > 0 $ is a given constant,
\mbox{$ \mbox{\boldmath $u$} = \mbox{\boldmath $u$}(x,t) $}
and
$ p = p(x,t) $
are the unknowns
(the flow velocity and pressure, respectively),
with condition (1.1$c$)
satisfied in
$ \mbox{\boldmath $L$}^{2}(\mathbb{R}^{n}) $,
i.e.,
$ {\displaystyle
\,\!\,\!
\|\, \mbox{\boldmath $u$}(\cdot,t) - \mbox{\boldmath $u$}_0 \;\!
\|_{L^{2}(\mathbb{R}^{n})}
\!\;\!\rightarrow 0
\,\!\,\!\,\!
} $
as $ \,\!\,\! t \,\mbox{\footnotesize $\searrow$}\,0 $.
In the present work,
we always assume $ \:\!2 \leq n \leq 4 $. \\

\mbox{} \vspace{-1.450cm} \\

A well known property of
Leray solutions
is that they are eventually
very regular:
there is always some
$ \!\;\!\;\!t_{\ast}\!\;\! \geq 0 \!\;\!\;\!$
such that one has
$ {\displaystyle
\!\;\!\;\!
\mbox{\boldmath $u$} \in
C^{\infty}(\:\!\mathbb{R}^{n} \!\times (\:\!t_{\ast}, \infty)\,\!)
} $
and, moreover, \\
\mbox{} \vspace{-0.600cm} \\
\begin{equation}
\tag{1.3}
\mbox{\boldmath $u$}(\cdot,t)
\in
C(\:\!(\:\!t_{\ast}, \infty), \:\!
\mbox{\boldmath $H$}\mbox{}^{\!\:\!m}\!\;\!(\mathbb{R}^{n})),
\quad \;\;\,
\forall \;\;\! m \geq 0,
\end{equation}
\mbox{} \vspace{-0.225cm} \\
see e.g.\;\cite{Kato1984, KreissLorenz1989, %
Leray1933, Leray1934, SchonbekWiegner1996, %
Sohr2001}.\footnote{%
%
%
%
%
It is known that
$ \;\!t_{\ast} \!\!\;\!\;\!= 0 \!\;\!\;\!$
if $ n = 2 $,
$ {\displaystyle
\;\!
t_{\ast}
\!\;\!\leq\!\;\!\;\!
\nu^{-\,5} \;\!
\|\, \mbox{\boldmath $u$}_{0} \;\!
\|_{\mbox{}_{\scriptstyle L^{2}(\mathbb{R}^{3})}}
  ^{\:\!4}
\!\!\!\;\!\;\!
} $
if $ n = 3 $,
$ {\displaystyle
\;\!
t_{\ast}
\!\;\!\leq\!\;\!\;\!
\nu^{-\,3} \;\!
\|\, \mbox{\boldmath $u$}_{0} \;\!
\|_{\mbox{}_{\scriptstyle L^{2}(\mathbb{R}^{4})}}
  ^{\:\!2}
\!\!\!\;\!\;\!
} $
if $ n = 4 $.
}
%
%
%
It is also well established that
%
%
$ {\displaystyle
\;\!
\lim_{t\,\rightarrow\,0}
\!\;\!\;\!
\|\, \mbox{\boldmath $u$}(\cdot,t) \,
\|_{\mbox{}_{\scriptstyle L^{2}(\mathbb{R}^{n})}}
\!\!\!\;\!\;\!=\!\;\!\;\! 0
\!\;\!\;\!
} $
and,
more generally, \\
\mbox{} \vspace{-0.850cm} \\
\begin{equation}
\tag{1.4}
\lim_{t\,\rightarrow\,0}
\,
t^{\:\!m/2}
\,
\|\, D^{m} \mbox{\boldmath $u$}(\cdot,t) \,
\|_{\mbox{}_{\scriptstyle L^{2}(\mathbb{R}^{n})}}
\!=\, 0
\end{equation}
\mbox{} \vspace{-0.200cm} \\
for every $ \:\!m \geq 1 $,
and for all Leray solutions
to the system (1.1)
\cite{BenameurSelmi2012, Kato1984, Masuda1984,
OliverTiti2000, Schonbek1985, SchonbekWiegner1996}.
Furthermore,
suitable stronger assumptions on the initial data
have led to interesting finer estimates
for the solutions and their derivatives,
see e.g.\;\cite{KreissHagstromLorenzZingano2003,
OliverTiti2000, SchonbekWiegner1996, Wiegner1987}.
An~important shortcut
for many of these results
(including (1.4) and the
Schonbek-Wiegner estimates
\cite{SchonbekWiegner1996})
is provided by
the following
fundamental inequality
recently discovered by the authors,
which has eluded previous studies.\footnote{%
%
%
%
%
For the definition of
$ {\displaystyle
\:\!
\|\, \mbox{\boldmath $u$}(\cdot,t) \,
\|_{\mbox{}_{\scriptstyle L^{2}(\mathbb{R}^{n})}}
\!\!\;\!\;\!
} $,
$ {\displaystyle
\|\, D^{m} \,\!\mbox{\boldmath $u$}(\cdot,t) \,
\|_{\mbox{}_{\scriptstyle L^{2}(\mathbb{R}^{n})}}
\!\:\!
} $
and other similar norms,
see (1.6).
}
%
%
\nl
\mbox{} \vspace{-0.250cm} \\
\mbox{} \hspace{-1.000cm}
\fbox{%
\mbox{} \hspace{+0.500cm}
\begin{minipage}[t]{14.900cm}
\nl
\mbox{} \vspace{-0.175cm} \\
%
%
%
%
{\bf Main Theorem.}
\textit{%
\,Let $ \,n \leq 4 $,
$ \mbox{\boldmath $u$}_{0} \in
\mbox{\boldmath $L$}^{2}_{\sigma}(\mathbb{R}^{n}) $,
and let $ \;\!\mbox{\boldmath $u$}(\cdot,t)$
be any particular
Leray solution to the Navier-Stokes equations $(1.1)$.
Then we have,
for every
$ \,\alpha \geq 0\:\!$\mbox{\em :} \\
}
\mbox{} \vspace{-0.525cm} \\
\begin{equation}
\notag
\mbox{} \;
\limsup_{t\,\rightarrow\,\infty}
\;
t^{\;\!\alpha\;\!+\;\!m/2} \,
\|\, D^{m} \mbox{\boldmath $u$}(\cdot,t) \,
\|_{\mbox{}_{\scriptstyle L^{2}(\mathbb{R}^{n})}}
\;\!\leq\,
K\!\;\!(\alpha, m) \:
\nu^{-\,m/2} \,
\limsup_{t\,\rightarrow\,\infty}
\;
t^{\:\!\alpha} \;\!
\|\, \mbox{\boldmath $u$}(\cdot,t) \,
\|_{\mbox{}_{\scriptstyle L^{2}(\mathbb{R}^{n})}}
\end{equation}
\mbox{} \vspace{-0.650cm} \\
\mbox{} \hfill (1.5) \\
\mbox{} \vspace{-0.575cm} \\
\textit{%
for every
$\, m \geq 1 $,
where
} \\
\mbox{} \vspace{-0.750cm} \\
\begin{equation}
\notag
\,
K\!\;\!(\alpha,m)
\,=\;
\mbox{$ {\displaystyle \min_{\delta \,>\,0} }$}
\:
\big\{\;\!
\delta^{\mbox{}^{\scriptstyle \:\!-\,1/2}}
\,\!
\prod_{\:\!j\,=\,0}^{\;\!m} \:\!
\bigl(\:\! \alpha + j/2 + \delta \;\!\bigr)^{\!\:\!1/2}
\:\!\bigr\}.
\end{equation}
\mbox{} \vspace{-0.400cm} \\
\end{minipage}
\mbox{} \hspace{+0.400cm} \mbox{} \\
}
%

%
%

\newpage

\mbox{} \vspace{-1.500cm} \\

In Section~2 we present
our original derivation
of (1.5),
which was based in part
on some previous ideas in
\cite{BrazSchutzZingano2013, KreissHagstromLorenzZingano2002, %
KreissHagstromLorenzZingano2003, Zingano1999}.
Alternative proofs
could also be developed \linebreak
(using e.g.\;Schonbek's Fourier splitting method
\cite{Schonbek1985, Schonbek1995}),
but we prefer
to follow \linebreak
the very way
in which (1.5)
was first revealed. \\
\nl
%
%
%
%
{\bf Notation.}
As already shown,
boldface letters
are used for
vector quantities,
as in
\linebreak
$ {\displaystyle
\mbox{\boldmath $u$}(x,t)
=
} $
$ {\displaystyle
(\:\! u_{\mbox{}_{\!\:\!1}}\!\;\!(x,t), ...\,\!\,\!,
 \:\! u_{\mbox{}_{\scriptstyle \!\;\! n}}\!\;\!(x,t) \:\!)
} $.
$\!$Also,
$ \nabla p \;\!\equiv \nabla p(\cdot,t) $
denotes the spatial gradient of $ \;\!p(\cdot,t) $;
$ \:\!D_{\!\;\!j} \!\;\!=\:\! \partial / \partial x_{\!\;\!j} \:\!$;
$ {\displaystyle
\:\!
\nabla \!\cdot \mbox{\boldmath $u$}
\:\!=
  D_{\mbox{}_{\!\:\!1}} u_{\mbox{}_{\!\:\!1}} \!\;\!+
  ... \!\;\!+
  D_{\mbox{}_{\scriptstyle \!\;\!n}} \,\!
  u_{\mbox{}_{\scriptstyle \!\;\!n}}
} $
is the (spatial) divergence of
$ \:\!\mbox{\boldmath $u$}(\cdot,t) $.
\linebreak
$ \!\;\!\mbox{\boldmath $L$}^{2}_{\sigma}(\mathbb{R}^{n}) $
denotes the
space
of solenoidal fields
$ \:\!\mbox{\bf v} = (v_{1}, ..., v_{n}) \!\:\!\in
\mbox{\boldmath $L$}^{2}(\mathbb{R}^{n}) $
$ \equiv L^{2}(\mathbb{R}^{n})^{n} \!\:\!$
with
$ \nabla \!\cdot \mbox{\bf v} \!\;\!= 0 $
in the distributional sense;
$ \dot{\mbox{\boldmath $H$}}\mbox{}^{1}(\mathbb{R}^{n}) =
\dot{H}^{1}(\mathbb{R}^{n})^{n} $
with
$ \dot{H}^{1}(\mathbb{R}^{n}) $
being
\linebreak
the homogeneous $L^{2}\!\!\;\!\;\!$
Sobolev space of order~1;
$ \mbox{\boldmath $H$}^{m}(\mathbb{R}^{n}) =
H^{m}(\mathbb{R}^{n})^{n} \!\:\!$,
where $ H^{m}(\mathbb{R}^{n}) $ is
the space of functions
$ v \in L^{2}(\mathbb{R}^{n}) $
whose $m$-th order derivatives
are also square~integrable.
$ C_{\mbox{\scriptsize w}}
(I, \:\!\mbox{\boldmath $L$}^{2}(\mathbb{R}^{n})) $
denotes the set of mappings from a
given interval~\mbox{$ I \subseteq \mathbb{R} $}
to
$ \mbox{\boldmath $L$}^{2}(\mathbb{R}^{n}) $
that are $L^{2}$-\;\!weakly continuous
at each $\:\! t \in I$.
$ {\displaystyle
\| \:\!\cdot\:\!
\|_{\scriptstyle L^{q}(\mathbb{R}^{n})}
\!\;\!
} $,
$ 1 \leq q \leq \infty $,
are the standard norms
of the Lebesgue spaces
$ L^{q}(\mathbb{R}^{n}) $,
with the vector counterparts \\
\mbox{} \vspace{-0.650cm} \\
\begin{equation}
\tag{1.6$a$}
\|\, \mbox{\boldmath $u$}(\cdot,t) \,
\|_{\mbox{}_{\scriptstyle L^{q}(\mathbb{R}^{n})}}
\;\!=\;
\Bigl\{\,
\sum_{i\,=\,1}^{n} \int_{\mathbb{R}^{n}} \!
|\:u_{i}(x,t)\,|^{q} \,dx
\,\Bigr\}^{\!\!\:\!1/q}
\end{equation}
\mbox{} \vspace{-0.675cm} \\
\begin{equation}
\tag{1.6$b$}
\mbox{} \;\;
\|\, D^{m} \mbox{\boldmath $u$}(\cdot,t) \,
\|_{\mbox{}_{\scriptstyle L^{q}(\mathbb{R}^{n})}}
\;\!=\;
\Bigl\{\!\!
\sum_{\mbox{} \;\;i, \,j_{\mbox{}_{1}} \!,..., \,j_{\mbox{}_{m}} =\,1}^{n}
\!\;\! \int_{\mathbb{R}^{n}} \!
|\, D_{\!\;\!j_{\mbox{}_{1}}}
\!\!\!\;\!\cdot\!\,\!\cdot\!\,\!\cdot \!\:\!
D_{\!\;\!j_{\mbox{}_{m}}}
\!\:\! u_{i}(x,t)\,|^{q} \,dx
\,\Bigr\}^{\!\!\:\!1/q}
\end{equation}
\mbox{} \vspace{-0.025cm} \\
if $ 1 \leq q < \infty $\/;
if $\, q = \infty $,
then
$ {\displaystyle
\;\!
\|\, \mbox{\boldmath $u$}(\cdot,t) \,
\|_{\mbox{}_{\scriptstyle L^{\infty}(\mathbb{R}^{n})}}
\!=\;\!
\max \, \bigl\{\,
\|\,u_{i}(\cdot,t)\,
\|_{\mbox{}_{\scriptstyle L^{\infty}(\mathbb{R}^{n})}}
\!\!: \, 1 \leq i \leq n
\,\bigr\}
} $, \linebreak
\mbox{} \vspace{-0.530cm} \\
$ {\displaystyle
\|\, D \,\!\mbox{\boldmath $u$}(\cdot,t) \,
\|_{\mbox{}_{\scriptstyle L^{\infty}(\mathbb{R}^{n})}}
\!=\;\!
\max \, \bigl\{\,
\|\, D_{\!\;\!j} \;\! u_{i}(\cdot,t)\,
\|_{\mbox{}_{\scriptstyle L^{\infty}(\mathbb{R}^{n})}}
\!\!: \:\! 1 \leq i, \:\!j \leq n
\,\bigr\}
\!\;\!
} $,
$\!\;\!\;\!$and
so forth. \\
\nl
%

%
%

{\bf 2. Proof of (1.5)} \\
\mbox{} \vspace{-0.650cm} \\

The derivation of (1.5) below
takes advantage of the
regularity property~(1.3)
and proceeds by induction in $m$.
It combines standard techniques
(energy inequalities
and related interpolation estimates)
with well known properties
of Leray solutions
(namely,
that
$ {\displaystyle
\!\;\!\;\!
\|\, \mbox{\boldmath $u$}(\cdot,t) \,
\|_{\mbox{}_{\scriptstyle L^{2}(\mathbb{R}^{n})}}
\!\rightarrow 0
} $
(as $\:\!t \rightarrow \infty$),
or that \\
\mbox{} \vspace{-0.600cm} \\
\begin{equation}
\tag{2.1}
\lim_{t\,\rightarrow\,\infty}
\,
\|\, D\:\!\mbox{\boldmath $u$}(\cdot,t) \,
\|_{\mbox{}_{\scriptstyle L^{2}(\mathbb{R}^{n})}}
=\;0,
\end{equation}
\mbox{} \vspace{-0.175cm} \\
which are easy to obtain directly).
\!As the proofs
for $ n = 2, 3, 4 $
are entirely similar,
we will present the details
for one case only --- say, $ n = 4 $.
Let
then
$ {\displaystyle
\:\!
\mbox{\boldmath $u$}(\cdot,t)
} $
be any given Leray solution
to~(1.1),
in $ \mathbb{R}^{4}\!\:\!$,
such that
we have,
for some $ \alpha \geq 0 $, \\
\mbox{} \vspace{-0.575cm} \\
\begin{equation}
\tag{2.2}
\limsup_{t\,\rightarrow\,\infty}
\;
t^{\:\!\alpha}
\;\!
\|\, \mbox{\boldmath $u$}(\cdot,t) \,
\|_{\mbox{}_{\scriptstyle L^{2}(\mathbb{R}^{4})}}
=:\:
\lambda_{0}(\alpha) \,<\, \infty.
\end{equation}
\mbox{} \vspace{-0.200cm} \\
Let $ \:\! \delta > 0 $,
$ 0 < \epsilon < 2 $ be given,
and
let $ \:\! t_{\ast} \!\;\!$
be the solution's regularity time
as defined in (1.3).
Recalling the basic estimate \\
\mbox{} \vspace{-0.675cm} \\
\begin{equation}
\tag{2.3}
\|\, \mbox{u}\,
\|_{\mbox{}_{\scriptstyle L^{4}(\mathbb{R}^{4})}}
\;\!\leq\;
\|\, D \:\!\mbox{u}\,
\|_{\mbox{}_{\scriptstyle L^{2}(\mathbb{R}^{4})}}
\!\:\!,
\end{equation}
\mbox{} \vspace{-0.500cm} \\
from which we get \\
\mbox{} \vspace{-0.625cm} \\
\begin{equation}
\tag{2.4}
\|\, D^{\ell} \,\!\mbox{u}\,
\|_{\mbox{}_{\scriptstyle L^{4}(\mathbb{R}^{4})}}
\|\, D^{m - \ell} \,\!\mbox{u}\,
\|_{\mbox{}_{\scriptstyle L^{4}(\mathbb{R}^{4})}}
\!\;\!
\leq\;
\|\, D \:\!\mbox{u}\,
\|_{\mbox{}_{\scriptstyle L^{2}(\mathbb{R}^{4})}}
\|\, D^{m + 1} \,\!\mbox{u}\,
\|_{\mbox{}_{\scriptstyle L^{2}(\mathbb{R}^{4})}}
\end{equation}
\mbox{} \vspace{-0.225cm} \\
for arbitrary $ \;\! m \geq 0 $,
$ 0 \leq \ell \leq m  $,
\;\!we may proceed
along the lines of
\cite{BrazSchutzZingano2013, Zingano1999}
as follows.
Taking the dot product
of (1.1$a$)
with
$ \:\!(\:\!t - t_{0})^{2\:\!\alpha \;\!+\;\! \delta}
\:\! \mbox{\boldmath $u$}(x,t) $
and integrating the result on
$ \:\!\mathbb{R}^{4} \!\:\!\times\!\;\![\,t_0, \:\!t\;\!] $,
for $ \:\! t \geq t_0 > t_{\ast} $,
we obtain,
because of (1.1$b$), \\
\mbox{} \vspace{-0.100cm} \\
\mbox{} \hspace{-0.250cm}
$ {\displaystyle
(\:\!t - t_{0})^{2\:\!\alpha \;\!+\;\!\delta}
\,
\|\, \mbox{\boldmath $u$}(\cdot,t) \,
\|_{\mbox{}_{\scriptstyle L^{2}(\mathbb{R}^{4})}}
  ^{\:\!2}
\!\;\!+\:
2\, \nu \!\!\;\!
\int_{\scriptstyle t_0}^{\;\!t}
\!\!\:\!
(\tau - t_{0})^{2\:\!\alpha \;\!+\;\!\delta}
\,
\|\, D \:\!\mbox{\boldmath $u$}(\cdot,\tau) \,
\|_{\mbox{}_{\scriptstyle L^{2}(\mathbb{R}^{4})}}
  ^{\:\!2}
d\tau
} $ \\
\mbox{} \vspace{-0.100cm} \\
\mbox{} \hfill
$ {\displaystyle
=\;
(2\:\!\alpha + \delta)
\!\!\;\!
\int_{\scriptstyle t_0}^{\;\!t}
\!\!\;\!
(\tau - t_{0})^{2\:\!\alpha \;\!+\;\!\delta \;\!-\;\!1}
\,
\|\, \mbox{\boldmath $u$}(\cdot,\tau) \,
\|_{\mbox{}_{\scriptstyle L^{2}(\mathbb{R}^{4})}}
  ^{\:\!2}
d\tau
} $ \\
\mbox{} \vspace{+0.050cm} \\
for $ \:\! t \geq t_0 \!\geq t_{\ast} $.
This promptly gives,
by (2.2),
that \\
\mbox{} \vspace{-0.600cm} \\
\begin{equation}
\tag{2.5}
\int_{\scriptstyle t_0}^{\;\!t}
\!\!\:\!
(\tau - t_{0})^{2\:\!\alpha \;\!+\;\!\delta}
\;\!
\|\, D \:\!\mbox{\boldmath $u$}(\cdot,\tau) \,
\|_{\mbox{}_{\scriptstyle L^{2}(\mathbb{R}^{4})}}
  ^{\:\!2}
d\tau
\:\leq\;
\frac{1}{2\;\!\nu}
\;
\frac{2\:\!\alpha \;\!+\;\!\delta}{\delta}
\:
(\lambda_{0}(\alpha) + \epsilon)^{2}
\,
(\:\!t - t_{0})^{\delta}
\end{equation}
\mbox{} \vspace{-0.100cm} \\
for all $ \!\;\!\;\!t \!\;\!\geq t_{0} \!\;\!$
(choosing
$ \!\;\!\;\!t_0 \geq t_{\ast} \!\;\!$
sufficiently large).
Next,
for $ \:\!m = 1 $,
we similarly have \\
%
%
%
\mbox{} \vspace{-0.100cm} \\
\mbox{} \hspace{-0.250cm}
$ {\displaystyle
(\:\!t - t_{0})^{2\:\!\alpha \;\!+\;\!1 \;\!+\;\!\delta}
\,
\|\, D \:\!\mbox{\boldmath $u$}(\cdot,t) \,
\|_{\mbox{}_{\scriptstyle L^{2}(\mathbb{R}^{4})}}
  ^{\:\!2}
\!\;\!+\:
2\, \nu \!\!\;\!
\int_{\scriptstyle t_0}^{\;\!t}
\!\!\:\!
(\tau - t_{0})^{2\:\!\alpha \;\!+\;\! 1 \;\!+\;\!\delta}
\,
\|\, D^{2} \,\!\mbox{\boldmath $u$}(\cdot,\tau) \,
\|_{\mbox{}_{\scriptstyle L^{2}(\mathbb{R}^{4})}}
  ^{\:\!2}
d\tau
} $ \\
\mbox{} \vspace{+0.000cm} \\
\mbox{} \hspace{+1.750cm}
$ {\displaystyle
\leq\;
(2\:\!\alpha + 1 + \delta)
\!\!\;\!
\int_{\scriptstyle t_0}^{\;\!t}
\!\!\;\!
(\tau - t_{0})^{2\:\!\alpha \;\!+\;\!\delta}
\,
\|\, D \:\!\mbox{\boldmath $u$}(\cdot,\tau) \,
\|_{\mbox{}_{\scriptstyle L^{2}(\mathbb{R}^{4})}}
  ^{\:\!2}
d\tau
\;\;\!+
} $ \\
\mbox{} \vspace{-0.025cm} \\
\mbox{} \hfill
$ {\displaystyle
K_{\mbox{}_{\!\:\!1}}
\!\!\,\!
\int_{\scriptstyle t_0}^{\;\!t}
\!\!\;\!
(\tau - t_{0})^{2\:\!\alpha \;\!+\;\! 1 \;\!+\;\!\delta}
\,
\|\, \mbox{\boldmath $u$}(\cdot,\tau) \,
\|_{\mbox{}_{\scriptstyle L^{4}(\mathbb{R}^{4})}}
\|\, D \:\!\mbox{\boldmath $u$}(\cdot,\tau) \,
\|_{\mbox{}_{\scriptstyle L^{4}(\mathbb{R}^{4})}}
\|\, D^{2} \,\!\mbox{\boldmath $u$}(\cdot,\tau) \,
\|_{\mbox{}_{\scriptstyle L^{2}(\mathbb{R}^{4})}}
d\tau
} $ \\
\mbox{} \vspace{+0.100cm} \\
(where
$ K_{\mbox{}_{\!\;\!1}} \!\!= 8 \,\sqrt{\:\!2\,} $),
which gives,
by (2.3): \\
\mbox{} \vspace{-0.100cm} \\
\mbox{} \hspace{-0.250cm}
$ {\displaystyle
(\:\!t - t_{0})^{2\:\!\alpha \;\!+\;\!1 \;\!+\;\!\delta}
\,
\|\, D \:\!\mbox{\boldmath $u$}(\cdot,t) \,
\|_{\mbox{}_{\scriptstyle L^{2}(\mathbb{R}^{4})}}
  ^{\:\!2}
\!\;\!+\:
2\, \nu \!\!\;\!
\int_{\scriptstyle t_0}^{\;\!t}
\!\!\:\!
(\tau - t_{0})^{2\:\!\alpha \;\!+\;\! 1 \;\!+\;\!\delta}
\,
\|\, D^{2} \,\!\mbox{\boldmath $u$}(\cdot,\tau) \,
\|_{\mbox{}_{\scriptstyle L^{2}(\mathbb{R}^{4})}}
  ^{\:\!2}
d\tau
} $ \\
\mbox{} \vspace{+0.000cm} \\
\mbox{} \hspace{+4.425cm}
$ {\displaystyle
\leq\;
(2\:\!\alpha + 1 + \delta)
\!\!\;\!
\int_{\scriptstyle t_0}^{\;\!t}
\!\!\;\!
(\tau - t_{0})^{2\:\!\alpha \;\!+\;\!\delta}
\,
\|\, D \:\!\mbox{\boldmath $u$}(\cdot,\tau) \,
\|_{\mbox{}_{\scriptstyle L^{2}(\mathbb{R}^{4})}}
  ^{\:\!2}
d\tau
\;\;\!+
} $ \\
\mbox{} \vspace{-0.025cm} \\
\mbox{} \hfill
$ {\displaystyle
K_{\mbox{}_{\!\;\!1}}
\!\!\,\!
\int_{\scriptstyle t_0}^{\;\!t}
\!\!\;\!
(\tau - t_{0})^{2\:\!\alpha \;\!+\;\! 1 \;\!+\;\!\delta}
\,
\|\, D \:\!\mbox{\boldmath $u$}(\cdot,\tau) \,
\|_{\mbox{}_{\scriptstyle L^{2}(\mathbb{R}^{4})}}
\|\, D^{2} \,\!\mbox{\boldmath $u$}(\cdot,\tau) \,
\|_{\mbox{}_{\scriptstyle L^{2}(\mathbb{R}^{4})}}
  ^{\:\!2}
d\tau
} $ \\
\mbox{} \vspace{+0.150cm} \\
for $ \:\!t \geq t_0 $.
By (2.1) and (2.5),
we then get
(increasing $ \:\!t_0 $ if necessary): 

\mbox{} \vspace{-0.950cm} \\
\begin{equation}
\tag{2.6$a$}
(\:\!t - t_{0})^{2\:\!\alpha \;\!+\;\! 1}
\,
\|\, D \:\! \mbox{\boldmath $u$}(\cdot,t) \,
\|_{\mbox{}_{\scriptstyle L^{2}(\mathbb{R}^{4})}}
  ^{\:\!2}
\,\leq\;
\frac{1}{2\;\!\nu}
\;
(2\:\!\alpha + 1 + \delta)
\;
\frac{2\:\!\alpha \;\!+\;\!\delta}{\delta}
\:
(\lambda_{0}(\alpha) + \epsilon)^{2}
\end{equation}
\mbox{} \vspace{-0.500cm} \\
and \\
\mbox{} \vspace{-0.650cm} \\
\begin{equation}
\tag{2.6$b$}
\int_{\scriptstyle t_0}^{\;\!t}
\!\!\;\!
(\tau - t_{0})^{2\:\!\alpha \;\!+\;\! 1 \;\!+\;\!\delta}
\;\!
\|\, D^{2} \,\!\mbox{\boldmath $u$}(\cdot,\tau) \,
\|_{\mbox{}_{\scriptstyle L^{2}(\mathbb{R}^{4})}}
  ^{\:\!2}
d\tau
\:\leq\:
\frac{(2\:\!\alpha + 1 + \delta)\;\!(2\:\!\alpha + \delta)}
{\delta \cdot [\;\!(2 - \epsilon)\!\;\!\;\!\nu \,]^{2}}
\;
(\lambda_{0}(\alpha) + \epsilon)^{2}
\,
(\:\!t - t_{0})^{\delta}
\end{equation}
\mbox{} \vspace{-0.500cm} \\
for all $\:\! t \geq t_{0} $.
\!Proceeding in this way
($ m = 2, 3, ...$)
we  obtain
at the $m\:\!$th step
 \\
%
%
%
\mbox{} \vspace{-0.050cm} \\
\mbox{} \hspace{-0.250cm}
$ {\displaystyle
(\:\!t - t_{0})^{2\:\!\alpha \;\!+\;\!m \;\!+\;\!\delta}
\,
\|\, D^{m} \,\!\mbox{\boldmath $u$}(\cdot,t) \,
\|_{\mbox{}_{\scriptstyle L^{2}(\mathbb{R}^{4})}}
  ^{\:\!2}
\!\;\!+\:
2\, \nu \!\!\;\!
\int_{\scriptstyle t_0}^{\;\!t}
\!\!\:\!
(\tau - t_{0})^{2\:\!\alpha \;\!+\;\! m \;\!+\;\!\delta}
\,
\|\, D^{m + 1} \,\!\mbox{\boldmath $u$}(\cdot,\tau) \,
\|_{\mbox{}_{\scriptstyle L^{2}(\mathbb{R}^{4})}}
  ^{\:\!2}
d\tau
} $ \\
\mbox{} \vspace{+0.000cm} \\
\mbox{} \hspace{+2.000cm}
$ {\displaystyle
\leq\;
(2\:\!\alpha + m + \delta)
\!\!\;\!
\int_{\scriptstyle t_0}^{\;\!t}
\!\!\;\!
(\tau - t_{0})^{2\:\!\alpha \;\!+\;\! m \;\!-\;\! 1 \;\!+\;\! \delta}
\,
\|\, D^{m} \,\!\mbox{\boldmath $u$}(\cdot,\tau) \,
\|_{\mbox{}_{\scriptstyle L^{2}(\mathbb{R}^{4})}}
  ^{\:\!2}
d\tau
\;\;\!+
} $ \\
\mbox{} \vspace{-0.225cm} \\
\mbox{} \hfill
$ {\displaystyle
\mbox{} \!
K_{\mbox{}_{\scriptstyle m}}
\!\!\,\!
\int_{\scriptstyle t_0}^{\;\!t}
\!\!\;\!
(\tau - t_{0})^{2\:\!\alpha \;\!+\;\! m \;\!+\;\!\delta}
\,
\|\, D^{m + 1} \,\!\mbox{\boldmath $u$}(\cdot,\tau) \,
\|_{\mbox{}_{\scriptstyle L^{2}(\mathbb{R}^{4})}}
\!\!\;\!
\sum_{\ell\,=\,0}^{[\,m/2\,]}
\!\!\;\!
\|\, D^{\ell} \,\! \mbox{\boldmath $u$}(\cdot,\tau) \,
\|_{\mbox{}_{\scriptstyle L^{4}(\mathbb{R}^{4})}}
\,\!
\|\, D^{m - \ell} \:\!\mbox{\boldmath $u$}(\cdot,\tau) \,
\|_{\mbox{}_{\scriptstyle L^{4}}}
\!\;\!\;\!
d\tau
} $ \\
\mbox{} \vspace{+0.125cm} \\
for $ t \geq t_{0} $,
and some constant
$ K_{\mbox{}_{\scriptstyle m}} \!\:\!> 0 $,
where
$ [\,m/2\,] $ denotes the
\mbox{integer part of $ \:\!m/2 $}.
This gives,
by (2.4): \\
\mbox{} \vspace{-0.125cm} \\
\mbox{} \hspace{-0.250cm}
$ {\displaystyle
(\:\!t - t_{0})^{2\:\!\alpha \;\!+\;\!m \;\!+\;\!\delta}
\,
\|\, D^{m} \,\!\mbox{\boldmath $u$}(\cdot,t) \,
\|_{\mbox{}_{\scriptstyle L^{2}(\mathbb{R}^{4})}}
  ^{\:\!2}
\!\;\!+\:
2\, \nu \!\!\;\!
\int_{\scriptstyle t_0}^{\;\!t}
\!\!\:\!
(\tau - t_{0})^{2\:\!\alpha \;\!+\;\! m \;\!+\;\!\delta}
\,
\|\, D^{m + 1} \,\!\mbox{\boldmath $u$}(\cdot,\tau) \,
\|_{\mbox{}_{\scriptstyle L^{2}(\mathbb{R}^{4})}}
  ^{\:\!2}
d\tau
} $ \\
\mbox{} \vspace{-0.300cm} \\
\mbox{} \hfill (2.7) \\
\mbox{} \vspace{-0.700cm} \\
\mbox{} \hspace{+2.000cm}
$ {\displaystyle
\leq\;
(2\:\!\alpha + m + \delta)
\!\!\;\!
\int_{\scriptstyle t_0}^{\;\!t}
\!\!\;\!
(\tau - t_{0})^{2\:\!\alpha \;\!+\;\! m \;\!-\;\! 1 \;\!+\;\! \delta}
\,
\|\, D^{m} \,\!\mbox{\boldmath $u$}(\cdot,\tau) \,
\|_{\mbox{}_{\scriptstyle L^{2}(\mathbb{R}^{4})}}
  ^{\:\!2}
d\tau
\;\;\!+
} $ \\
\mbox{} \vspace{-0.025cm} \\
\mbox{} \hfill
$ {\displaystyle
\mbox{} \!
\Bigl(\;\! 1 + \Bigl[\,
\mbox{\small $ {\displaystyle \frac{m}{2} }$} \,\Bigr]
\;\!\Bigr)
\cdot
K_{\mbox{}_{\scriptstyle m}}
\!\!\,\!
\int_{\scriptstyle t_0}^{\;\!t}
\!\!\;\!
(\tau - t_{0})^{2\:\!\alpha \;\!+\;\! m \;\!+\;\!\delta}
\,
\|\, D \:\!\mbox{\boldmath $u$}(\cdot,\tau) \,
\|_{\mbox{}_{\scriptstyle L^{2}(\mathbb{R}^{4})}}
\,\!
\|\, D^{m + 1} \,\!\mbox{\boldmath $u$}(\cdot,\tau) \,
\|_{\mbox{}_{\scriptstyle L^{2}(\mathbb{R}^{4})}}
  ^{\:\!2}
\!\;\!\;\!
d\tau
} $. \\
\mbox{} \vspace{+0.250cm} \\
At this stage,
we  would already know
from the previous steps
that \\
\mbox{} \vspace{-0.750cm} \\
\begin{equation}
\tag{2.8$a$}
(\!\;\!\;\!t - t_{0})^{2\!\;\!\;\!\alpha \;\!+\;\! k}
\;\!
\|\, D^{k} \,\! \mbox{\boldmath $u$}(\cdot,t) \,
\|_{\mbox{}_{\scriptstyle L^{2}(\mathbb{R}^{4})}}
  ^{\:\!2}
\!\;\!\leq\;
\frac{1}{\delta \cdot [\;\!(2 - \epsilon)\!\;\!\;\!\nu \,]^{k}}
\;\!\;\!
\biggl\{\;\!
\prod_{j\,=\,0}^{k}
(2\!\;\!\;\!\alpha + j + \delta)
\!\;\!\;\!
\biggr\}
\;\!
(\lambda_{0}(\alpha) + \epsilon)^{2}
\end{equation}
\mbox{} \vspace{-0.775cm} \\
and \\
\mbox{} \vspace{-0.120cm} \\
\mbox{} \hspace{-0.050cm}
$ {\displaystyle
\int_{\scriptstyle t_0}^{\;\!t}
\!\!\;\!
(\tau - t_{0})^{2\!\;\!\;\!\alpha \;\!+\;\! k \;\!+\;\!\delta}
\;\!
\|\, D^{k + 1} \mbox{\boldmath $u$}(\cdot,\tau) \,
\|_{\mbox{}_{\scriptstyle L^{2}(\mathbb{R}^{4})}}
  ^{\:\!2}
d\tau
\,\leq\;\!\;\!
\frac{\delta^{-\,1}}{\, [\;\!(2 - \epsilon)\!\;\!\;\!\nu \,]^{k + 1}}
\;\!\;\!
\biggl\{\;\!
\prod_{j\,=\,0}^{k}
(2\!\;\!\;\!\alpha + j + \delta)
\!\;\!\;\!
\biggr\}
} $ \\
\mbox{} \vspace{-0.200cm} \\
\mbox{} \hfill (2.8$b$) \\
\mbox{} \vspace{-0.675cm} \\
\mbox{} \hspace{+8.700cm}
$ {\displaystyle
\times \;
(\lambda_{0}(\alpha) + \epsilon)^{2}
\!\cdot
(\:\!t - t_{0})^{\delta}
} $ \\
\mbox{} \vspace{-0.025cm} \\
for all $ t \geq t_{0} $,
and each $ 0 \leq k < m $.
\!By (2.1) and (2.7),
and increasing $ t_0 $ if necessary,
we would then obtain
(2.8) for $ \:\!k = m \:\!$ as well,
completing the induction step. \linebreak

\mbox{} \vspace{-1.150cm} \\

The argument above established that,
for each $ \:\! m \geq 1 $,
we have \\
\mbox{} \vspace{-0.725cm} \\
\begin{equation}
\notag
(\!\;\!\;\!t - t_{0})^{2\!\;\!\;\!\alpha \;\!+\;\! m}
\;\!
\|\, D^{m} \,\! \mbox{\boldmath $u$}(\cdot,t) \,
\|_{\mbox{}_{\scriptstyle L^{2}(\mathbb{R}^{4})}}
  ^{\:\!2}
\!\;\!\leq\;
\frac{1}{\delta \cdot [\;\!(2 - \epsilon)\!\;\!\;\!\nu \,]^{m}}
\;\!\;\!
\biggl\{\;\!
\prod_{j\,=\,0}^{m}
(2\!\;\!\;\!\alpha + j + \delta)
\!\;\!\;\!
\biggr\}
\;\!
(\lambda_{0}(\alpha) + \epsilon)^{2}
\end{equation}
\mbox{} \vspace{-0.125cm} \\
for all $ \:\!t $
sufficiently large.
Since $ \:\! \delta > 0 $,
$ 0 < \epsilon < 2 \:\!$
are arbitrary,
this gives the result. \\
%

%
%
%
\mbox{} \vspace{-0.400cm} \\
%
%
%
%
{\bf Acknowledgements.} \\
Work of the first author was supported by
NSF Grant DMS-1418871,
and that of the last author by
CAPES Grant \mbox{\small \#}\,88881.067966/2014-01.
Any opinions, findings, and conclusions or recommendations
are those of the authors and do not necessarily reflect
the views of the National Science Foundation
or CAPES. \\
\mbox{} \vspace{-1.000cm} \\
%

%
%
%

%
%


\begin{thebibliography}{999}

%
%


\bibitem{BenameurSelmi2012}
{\sc J. Benameur and R. Selmi},
{\em Long time decay to the Leray solution of the
two-dimensional Navier-Stokes equations},
Bull. London Math. Soc. {\bf 44} (2012),
1001-1019.

\bibitem{BrazSchutzZingano2013}
{\sc P.\;Braz e Silva, L.\;Sch\"utz and P.\;R.\;Zingano},
{\it On some energy inequal-\linebreak
ities and supnorm estimates
for advection-diffusion equations in $\mathbb{R}^{n} \!\:\!$},
Nonl. Anal. {\bf 93} (2013),
90-96.

\bibitem{FujitaKato1964}
{\sc H. Fujita and T. Kato},
{\it On the Navier-Stokes initial value problem},
Arch. Rat. Mech. Anal. {\bf 16} (1964),
269-315.



\bibitem{Kato1984}
{\sc T. Kato},
{\it Strong $L^p$-$\:\!$solutions of the Navier\;\!-\:\!Stokes equations
in $\mathbb{R}^{m}\!\,\!$, with applications to weak solutions},
Math. Z. {\bf 187} (1984), 471-480.

\bibitem{KreissHagstromLorenzZingano2002}
{\sc H.-O. Kreiss, T. Hagstrom, J. Lorenz and P. R. Zingano},
{\it Decay in time of the solutions of the Navier-Stokes equations
for incompressible flows},
unpublished note,
University of New Mexico, Albuquerque, NM, 2002.

\bibitem{KreissHagstromLorenzZingano2003}
{\sc H.-O. Kreiss, T. Hagstrom, J. Lorenz and P. R. Zingano},
{\it Decay in time of incompressible flows},
J. Math. Fluid Mech. {\bf 5} (2003),
231-244.

\bibitem{KreissLorenz1989}
{\sc H.-O. Kreiss and J. Lorenz},
{\sf Initial-boundary value problems and the Navier-Stokes equations},
Academic Press,
New York, 1989.
(Reprinted in the series SIAM Classics in Applied Mathematics,
 Vol.\;47, 2004.)

\bibitem{Ladyzhenskaya1969}
{\sc O. A. Ladyzhenskaya},
{\sf The Mathematical Theory of Viscous Incompressible Flow} (2nd ed.),
Gordon and Breach, New York, 1969.

\bibitem{Leray1933}
{\sc J. Leray},
{\it \'Etude de diverses \'equations int\'egrales non lin\'eaires
et de quelques probl\`emes que pose l'hydrodynamique},
J. Math. Pures Appl. {\bf 12} (1933),
1-82.

\bibitem{Leray1934}
{\sc J. Leray},
{\it Essai sur le mouvement d'un fluide visqueux emplissant l'espace},
Acta Math. {\bf 63} (1934), 193-248.



\bibitem{Masuda1984}
{\sc K.$\;$Masuda},
{\it Weak solutions of the Navier-Stokes equations},
T\^ohoku Math. Journal {\bf 36} (1984),
623-646.

\bibitem{OliverTiti2000}
\textsc{M. Oliver and E. S. Titi},
{\it Remark on the rate of decay of higher order derivatives
for solutions to the Navier-Stokes equations in $ \mathbb{R}^{n} \!$},
J. Funct. Anal. {\bf 172} (2000),
1-18.

\bibitem{Schonbek1985}
{\sc M. E. Schonbek},
{\it $L^{2}$ decay for weak solutions
of the Navier-Stokes equations},
Arch. Rat. Mech. Anal. {\bf 88} (1985),
209-222.


\bibitem{Schonbek1995}
{\sc M. E. Schonbek},
{\it The Fourier splitting method},
in:
P. Concus and K. Lancaster (Eds.),
{\sf Advances in Geometric Analysis and Continuum Mechanics},
International Press, Cambridge, 1995, pp.$\;$269-274.


\bibitem{SchonbekWiegner1996}
{\sc M. E. Schonbek and M. Wiegner},
{\it On the decay of higher-order norms of the
solutions of Navier-Stokes equations},
Proc. Roy. Soc. Edinburgh {\bf 126A} (1996),
677-685.


\bibitem{ShinbrotKaniel1966}
{\sc M. Shinbrot and S. Kaniel},
{\it The initial value problem
for the Navier-Stokes equations},
Arch. Rat. Mech. Anal. {\bf 21} (1966),
270-285.

\bibitem{Sohr2001}
{\sc H. Sohr},
{\sf The Navier-Stokes Equations},
Birkh\"auser, Basel, 2001.

\bibitem{Temam1984}
{\sc R. Temam},
{\sf Navier-Stokes Equations:
theory and numerical analysis} (2nd ed.),
AMS/Chelsea,
Providence, 1984.


\bibitem{Wiegner1987}
{\sc M. Wiegner},
{\it Decay results for weak solutions of the Navier-Stokes
equations on $ \mathbb{R}^{n}\!$},
J. London Math. Soc. {\bf 35} (1987), 303-313.


\bibitem{Zingano1999}
{\sc P. R. Zingano},
{\it Nonlinear $L^{2}\!$ stability under large disturbances},
J. Comp. Appl. Math. {\bf 103} (1999),
207-219.

\end{thebibliography}
\end{document}